\documentclass[11pt,leqno,twoside]{article}

\usepackage{amsfonts,amsmath,amsthm,amssymb}
\usepackage{mathtools}
\usepackage{enumerate}
\usepackage{graphics,graphicx,subfigure}

\usepackage{color}
\usepackage{todonotes}
\usepackage{cancel}
\usepackage{url}
\usepackage{hyperref}
\usepackage{makeidx}
\usepackage{showidx}
\usepackage{multicol}        
\usepackage{xspace}
\usepackage{stmaryrd}        
\usepackage{pifont}          
\usepackage{fancybox}        
\usepackage{bm}
\usepackage{mathrsfs}

 \usepackage{fancyhdr}
\theoremstyle{plain}
 \theoremstyle{definition}
 \newtheorem{lem}{Lemma}
 \newtheorem{defn}[lem]{Definition}
 \newtheorem{thm}[lem]{Theorem}
 \newtheorem{prop}[lem]{Proposition}
 \newtheorem{cor}[lem]{Corollary}
 \newtheorem{notn}[lem]{Notations}
 \newtheorem{pb}[lem]{Problem}
 \newtheorem{form}[lem]{Formulae}
 
 \newtheorem*{rk}{Remark}
 \newtheorem*{com}{Comment}
 \newtheorem*{ex}{Example}
 \theoremstyle{remark}

 \newcommand{\blem}{\begin{lem}}
 \newcommand{\elem}{\end{lem}}
 \newcommand{\bdefn}{\begin{defn}}
 \newcommand{\edefn}{\end{defn}}
 \newcommand{\bthm}{\begin{thm} }
 \newcommand{\ethm}{\end{thm}}
 \newcommand{\bprop}{\begin{prop}}
 \newcommand{\eprop}{\end{prop}}
 \newcommand{\bcor}{\begin{cor}}
 \newcommand{\ecor}{\end{cor}}
 \newcommand{\bnotn}{\begin{notn}}
 \newcommand{\enotn}{\end{notn}}
 \newcommand{\bpb}{\begin{pb}}
 \newcommand{\epb}{\end{pb}}
 \newcommand{\bform}{\begin{form}}
 \newcommand{\eform}{\end{form}}
 \newcommand{\brk}{\begin{rk}}
 \newcommand{\erk}{\end{rk}}
 \newcommand{\bcom}{\begin{com}}
 \newcommand{\ecom}{\end{com}}
 \newcommand{\bex}{\begin{ex}}
 \newcommand{\eex}{\end{ex}}
 \newcommand{\bpf}{\begin{proof}}
 \newcommand{\epf}{\end{proof}}

\newcommand{\cF}{\mathcal{F}}

\newcommand{\cM}{\mathcal{M}}

\newcommand{\cP}{\mathcal{P}}

\newcommand{\cS}{\mathcal{S}}

\newcommand{\cU}{\mathcal{U}}
\newcommand{\cV}{\mathcal{V}}

\newcommand{\bN}{\mathbb{N}}

\newcommand{\bR}{\mathbb{R}}

\newcommand{\be}{\begin{equation}}
\newcommand{\ee}{\end{equation}}
\newcommand{\bal}{\begin{align}}
\newcommand{\eal}{\end{align}}
\newcommand{\ba}{\begin{align*}}
\newcommand{\ea}{\end{align*}}
\newcommand{\bmx}{\begin{matrix}}
\newcommand{\emx}{\end{matrix}}
\newcommand{\bbmx}{\begin{bmatrix}}
\newcommand{\ebmx}{\end{bmatrix}}
\newcommand{\bpmx}{\begin{pmatrix}}
\newcommand{\epmx}{\end{pmatrix}}
\newcommand{\bvmx}{\begin{vmatrix}}
\newcommand{\evmx}{\end{vmatrix}}

\newcommand{\wt}{\widetilde}
\newcommand{\f}{\frac}

\newcommand{\inc}{\subseteq}

\newcommand{\sgn}{\mathrm{sgn}}

\newcommand{\minimize}[1]{\underset{#1}{\rm minimize}\,}

\newcommand{\la}{\lambda}

\newcommand{\eps}{\varepsilon}
  
\def\N{\mathbb{N}}

\newcommand{\rev}[1]{{#1}}
\newcommand{\revv}[1]{{#1}}

\title{\vspace{-20mm}Optimization-Aided Construction of \\  Multivariate Chebyshev Polynomials\medskip\hrule height 1.2pt \vspace{-6mm}}
\author{M.  Dressler\footnote{School of Mathematics and Statistics, University of New South Wales,  Australia, \url{m.dressler@unsw.edu.au}}, S. Foucart\footnote{Department of Mathematics, Texas A\&M University, United States, \url{foucart@tamu.edu}}, 
\rev{M.  Joldes}\footnote{LAAS-CNRS, France, \url{joldes@laas.fr}},
\rev{E. de Klerk}\footnote{Department of Econometrics and Operations Research, Tilburg University, The Netherlands, \url{E.deKlerk@tilburguniversity.edu}}, 
 J.~B. Lasserre\footnote{LAAS-CNRS and Toulouse School of Economics, France, \url{lasserre@laas.fr}}, and Y. Xu\footnote{Department of Mathematics, University of Oregon, United States, \url{yuan@uoregon.edu}}}
\date{\vspace{-8mm}\rule{100mm}{0.8pt}}

\newcommand\shorttitle{Optimization-Aided Construction of Multivariate Chebyshev Polynomials}
\newcommand\authors{M.  Dressler, S. Foucart, \rev{M.  Joldes, E. de Klerk}, J.~B. Lasserre, and Y. Xu}

\begin{document}
\maketitle

\begin{abstract}
This article is concerned with an extension of  univariate Chebyshev polynomials of the first kind to the multivariate setting,
where one chases best approximants to specific monomials by polynomials of lower degree relative to the uniform norm.
Exploiting the Moment-SOS hierarchy,
we devise a versatile semidefinite-programming-based procedure to compute such best approximants,
as well as associated signatures.
Applying this procedure in three variables
leads to the values of best approximation errors for all monomials up to degree six on the euclidean ball, the simplex, and the cross-polytope.
Furthermore,
inspired by numerical experiments,
we obtain explicit expressions for Chebyshev polynomials in two cases unresolved before,
namely for the monomial $x_1^2 x_2^2 x_3$ on the euclidean ball 
and for the monomial $x_1^2 x_2 x_3$ on the simplex.
\end{abstract}



\section{Introduction}

For $n \ge 1$,  let $\cP_n$ denote the space of univariate polynomials of degree less than or equal to $n$.
The classical $n$th degree Chebyshev polynomial 
(of the first kind)
is the polynomial $T_n$ often implicitly defined via the relation $T_n(\cos \theta) = \cos(n \theta)$ for all $\theta \in [-\pi,\pi]$.
It is characterized by a wealth of extremal properties, including:
\vspace{-5mm}
\begin{itemize}
\item $2^{-n+1} T_n$ is the monic polynomial that deviates least from zero in the uniform norm on $[-1,1]$, i.e.,
$T_n$ minimizes $\|p\|_{[-1,1]} \coloneqq \max \{ |p(x)| : x \in [-1,1] \}$ over all polynomials $p \in \cP_n$ satisfying ${\rm coeff}_{x^n}(p) = 2^{n-1}$---this is how Chebyshev polynomials were first introduced in~\cite{CheOri};
\item $2^{-n+1} T_n$ is the monic polynomial that deviates least from zero in the $L_2$-norm on $[-1,1]$ with respect to the inverse semicircle weight,
i.e.,  $T_n$ minimizes $\int_{-1}^1 p(x)^2 (1-x^2)^{-1/2} dx$ over all polynomials $p \in \cP_n$ satisfying ${\rm coeff}_{x^n}(p) = 2^{n-1}$---this relates to the orthogonality of Chebyshev polynomials for this weight;
\item $T_n$ is the extremizer of every differentiation operator,
i.e.,  $T_n$ maximizes $\|p^{(k)}\|_{[-1,1]}$ over all polynomials $p \in \cP_n$ satisfying $\|p\|_{[-1,1]} \le 1$ for every $k=1,2,\ldots,n$---this is Markov's inequality due to A. A. Markov
for $k=1$ and Markov's inequality due to his younger brother V. A. Markov for $k=2,\ldots,n$;
\item $T_n$ is the polynomial with the largest growth outside $[-1,1]$,
i.e., \rev{for every $t \not\in [-1,1]$ and every $k=0,1,\ldots,n$,
  $T_n$ maximizes $| p^{(k)} (t) | $ over all polynomials $p \in \cP_n$ satisfying $\|p\|_{[-1,1]} \le 1$};
\item $T_n$ is the polynomial with largest arc-length on $[-1,1]$,
i.e.,  $T_n$ maximizes $\int_{-1}^1  \sqrt{1+p'(x)^2} dx$ over all polynomials $p \in \cP_n$ satisfying $\|p\|_{[-1,1]} \le 1$. 
\end{itemize}

Each of these five properties,
which are all found in the classic book \cite{Riv} by Rivlin,
could serve as a rationale for a generalization of Chebyshev polynomials to the multivariate setting.
The generalization examined in this article is based on the first property.
Thus, denoting by $\cP_n^d$ the space of $d$-variate polynomials of degree $\le n$
and considering a domain $\Omega \inc \bR^d$,
we intend to tackle the optimization program
$$
\minimize{p \in \cP_n^d} \; \|p\|_\Omega \coloneqq \max_{x \in \Omega} |p(x)| 
\qquad \mbox{subject to } p \mbox{ being monic.}
$$
Although this is a convex optimization program---the constraint is linear and the objective function is convex---solving it is far from trivial.
By introducing a slack variable $c \in \bR$,
it is seen to be equivalent to
\be
\label{ReformNonNegCst}
\minimize{c \in \bR, \, p \in \cP_n^d} \; c
\qquad \mbox{subject to } \mbox{$p$ being monic}
\mbox{ and to } \left\{ \begin{matrix}  c + p \ge 0 & \mbox{ on } \Omega,\\
c - p \ge 0 & \mbox{ on } \Omega.  \end{matrix} \right.
\ee
The added constraints are polynomial nonnegativity constraints and, as such,
can conceivably be dealt with using sum-of-squares (SOS) techniques.
In fact, 
\rev{as clarified later},
we will advantageously use the dual facet of the Moment-SOS hierarchy \cite{JBL} to address our central optimization program.

But before delving into the technicalities,
let us mention that the above formulation comes with some ambiguities about\vspace{-5mm}
\begin{itemize}
\item the notion of multivariate degree: we will concentrate exclusively on the total degree given by
$$
\deg \Bigg( \sum_{k = (k_1,\ldots,k_d)} c_{k_1,\ldots,k_d} \, x_1^{k_1} \cdots x_d^{k_d} \Bigg) = \max_{\rev{k: c_k \not= 0}} |k|,
\qquad \mbox{where } |k| = k_1 + \cdots+k_d;
$$
\item the choice of domain $\Omega$:
we will consider only the simplex $S$,
the cross-polytope $C$ ($\ell_1$-ball),
the euclidean ball $B$ ($\ell_2$-ball),
and the hypercube $H$ ($\ell_\infty$-ball),
which are given by
\begin{align*}
& S = \Big\{ x \in \bR^d: x_1,\ldots,x_d \ge 0 \mbox{ and } \sum_{i=1}^{d} x_i  \le 1 \Big\},
& & \quad C =  \Big\{ x \in \bR^d: \sum_{i=1}^{d} |x_i|  \le 1 \Big\},\\
& B = \Big\{ x \in \bR^d: \sum_{i=1}^{d} |x_i|^2  \le 1 \Big\},
& & \quad H = \Big\{ x \in \bR^d: \max_{i=1,\ldots,d} |x_i|  \le 1 \Big\};
\end{align*}
\item the meaning of `monic' in the constraint:
it can be interpreted as imposing that the coefficient on a fixed $n$th degree monomial \rev{$m_k$} equals one
while the coefficients on all other $n$th degree monomials equal zero,
leading to the best approximation problem
\be
\label{Opt1}
\minimize{p \in \cP_{n-1}^d} \; \|m_k - p\|_\Omega,
\qquad \mbox{where} \quad
m_k(x) = x_1^{k_1} \cdots x_d^{k_d}
\quad \mbox{and} \quad |k|=n,
\ee 
or it could be interpreted as imposing that the coefficients on all the $n$th degree monomials sum up to one,  leading to the program
\be
\label{Opt2}
\minimize{p \in \cP_{n}^d} \; \|p\|_\Omega
\qquad \mbox{subject to } \sum_{|k| = n} {\rm coeff}_{m_k}(p) = 1.
\ee
The term Chebyshev polynomial will refer to the first interpretation.
It is the subject of this article and necessitates a computational approach.
The second interpretation comes,
more classically,
with explicit expressions for a large class of domains including $S$, $C$, $B$, and $H$.
This is the subject of a companion article, see \cite{DFKJLX}. 
\end{itemize}

Here, our contribution includes the numerical---sometimes explicit---construction of all Chebyshev polynomials
up to degree $n=6$ for $d=3$ variables.
The whole list of errors of best approximation  is assembled in Section~\ref{SecEx},
completing a partial catalog of known results recalled in Section~\ref{SecKno}.
This section also provides a refresher on some important reductions and on the central concept of signature.
The production of novel Chebyshev polynomials exploits a semidefinite programming procedure presented in Section~\ref{SecDes}.
Arguably,  this is the centerpiece of our work
and we emphasize its versatility, 
which would allow one to make easy adjustments for related problems,
e.g. multivariate Zolotarev's polynomials could be constructed with only small modifications of the procedure.
It is also worth noting already at this point that the workflow is not only numerical:
the experimental Chebyshev polynomials returned by our procedure can be verified explicitly or symbolically to be genuine Chebyshev polynomials.
For instance, best approximants to the monomial $m_{(2,2,1)}$ relative to the euclidean ball and to the monomial $m_{(2,1,1)}$ relative to the simplex are derived analytically in Section \ref{SecEx}.
Finally, Section \ref{SecConc} gives an outlook on a possible augmentation of the procedure and its deployment into further computational endeavors.

\section{Prior \rev{Knowledge}}
\label{SecKno}

This section is exclusively concerned with (monomial-specific) multivariate Chebyshev polynomials,
i.e., with solutions to the best approximation problem \eqref{Opt1}.
We start by recalling a characterization of these solutions  involving the notion of signature.
Then we provide a catalog of previously derived multivariate Chebyshev polynomials---more precisely,  of the known results that we are aware of.
We point out from the outset that multivariate Chebyshev polynomials are generically not unique,
explaining our tendency to manipulate signatures preferably to polynomials themselves.

\subsection{Characterization via signatures}

To be most general,
let us assume that $\Omega$ is a compact set and that we are trying to approximate a continuous function $f \in C(\Omega)$ by elements of a finite-dimensional vector space $\cV \inc C(\Omega)$,
assuming of course that $f \not\in \cV$.
We have in mind the case where $f$ is a $n$th degree monomial and where $\cV$ is the space $\cP_{n-1}^d$,
but the considerations of this subsection are valid beyond this specific case.
The error of best approximation and (any one of) the best approximant(s) shall be denoted by 
$E_{\cV}(f,\Omega)$ and by $v^*$, respectively, so that
\be
\label{PrimalForm}
E_{\cV}(f,\Omega) = \min_{v \in \cV} \|f-v\|_\Omega = \|f-v^*\|_\Omega.
\ee
\rev{A folklore result in Approximation Theory shows that the latter can alternatively be expressed as a maximum, namely as}
\be
\label{DualForm}
E_{\cV}(f,\Omega)
= \max_{\la \in C(\Omega)^*} \; \la(f)
\qquad \mbox{subject to } \la_{| \cV} = 0 \mbox{ and } \|\la\|_{C(\Omega)^*} = 1.
\ee
\rev{In optimization custom, this is viewed as a strong duality result,
whose proof is informative to sketch here.
To start, the inequality that the maximum in \eqref{DualForm} does not exceed the minimum in \eqref{PrimalForm} simply follows from the fact that,
for any feasible $v \in \cV$ and $\la \in C(\Omega)^*$,
one has
$$
\la(f) = \la(f-v) \le \|\la\|_{C(\Omega)^*} \|f-v\|_\Omega = \|f-v\|_\Omega.  
$$
Next, for the reverse inequality,
given a best approximant $v^*$,
consider the linear functional $\wt{\la}$ defined on $\cV_f \coloneqq \cV \oplus \bR f$ by 
$$
\wt{\la}(v+tf) = t \|f-v^*\|_\Omega,
\quad v \in \cV, \quad t \in \bR.
$$
This linear functional has norm one:
indeed, setting aside the trivial case $t=0$, one has
$$
|\wt{\la}(v+tf)| = |t| \, \|f-v^*\|_\Omega
\le |t| \, \|f - (-v/t)\|_\Omega =  \|v+ tf \|_\Omega,
$$
with equality for $t=1$ and $v = - v^*$.
Therefore, by the Hahn--Banach extension theorem,
there exists a linear functional $\la$ defined on $C(\Omega)$ such that $\| \la\|_{C(\Omega)^*} = 1$ and $\la_{| \cV_f} = \wt{\la}_{| \cV_f}$, implying in particular $\la_{| \cV} = 0$.
The maximum in \eqref{DualForm} is then larger than or equal to $\la(f) = \wt{\la}(f)~=~\|f - v^*\|_\Omega$,
i.e., than the minimum in \eqref{PrimalForm}.
This concludes the sketch of strong duality. 
}

\rev{It is worth commenting on the form of the Hahn--Banch extension $\la$ in the above argument.
To do~so, we notice that the linear functional $\wt{\la} \in \cV_f^*$ is expressed as a convex combination of $L$ extreme points of the unit ball of $\cV_f^*$,
where one can take $L \le \dim(\cV_f)+1$ by applying Krein--Milman and Carath\'eodory theorems.
Since the extreme points of the unit ball of $\cV_f^*$ are restrictions to $\cV_f$ of some
$\pm \delta_\omega$, $\omega \in \Omega$,  where $\delta_\omega$ represents the Dirac evaluation functional at a point $\omega$,
we can write $\wt{\la} = \sum_{\ell =1}^L \tau_\ell \eps_\ell {\delta_{\omega_\ell}}_{|\cV_f}$
with $\omega_1,\ldots,\omega_L \in \Omega$,
$\eps_1, \ldots, \eps_L = \pm 1$,
and $\tau_1,\ldots, \tau_L > 0$ satisfying $\sum_{\ell=1}^L \tau_\ell = 1$. 
It is then clear that the Hahn--Banach extension $\la$ can be chosen as $\la = \sum_{\ell =1}^L \tau_\ell \eps_\ell \delta_{\omega_\ell}$.
Of note, the estimation of $L$ can be refined to $L \le \dim(\cV_f)$
by calling upon Theorem 2.13 from Rivlin's book \cite{Riv},
which states that a norm-one linear functional on a real finite-dimensional subspace $\cU$ of $C(\Omega)$ can be expressed as a convex combination of $L \le \dim(\cU)$ extreme points of the unit ball of $\cU^*$.
In the particular case $\cU = \cP_k^d$, 
it equates to Tchakaloff's theorem and its generalizations (notably by Richter \cite{Rich}), stating that if a measure $\mu$ on $\bR^d$ has moments up to degree $k$, then there exists an atomic measure with at most ${k+d \choose d}$ atoms in $\mathrm{supp}(\mu)$ and with same moments up to degree $k$.
}

\rev{Closing the digression on the value of $L$ and
keeping the above notation,
we point out that }
\begin{align}
\label{IntroSigna}
E_\cV(f,\Omega) & = \la(f) = \la(f-v^*)
= \sum_{\ell = 1}^L \tau_\ell \eps_\ell (f-v^*)(\omega_\ell)\\
\nonumber
& \le \sum_{\ell = 1}^L \tau_\ell |  (f-v^*)(\omega_\ell) | 
\le \sum_{\ell = 1}^L \tau_\ell \|f-v^*\|_\Omega 
= \|f-v^*\|_\Omega\\ 
\nonumber
& = E_\cV(f,\Omega),
\end{align}
 which implies that the equalities $\eps_\ell (f-v^*)(\omega_\ell) = \|f-v^*\|_\Omega  = E_\cV(f,\Omega)$ hold for all $\ell = 1,\ldots,L$.
 This brings us to the notion of extremal signature associated with $f-v^*$,
 defined below along the lines of \cite[Section 2.2]{Riv}.
 
\bdefn
\label{def:signature}
A signature with support (aka base) $\cS \inc \Omega$ is simply a function from $\cS$ to $\{\pm 1\}$.
A~signature $\sigma$ with support $\cS$ is said to be extremal for $\cV$ if there exist weights $\tau_\omega > 0$, $\omega \in \cS$, such that $\sum_{\omega \in \cS} \tau_\omega \sigma(\omega) v(\omega) = 0$ for all $v \in \cV$.
A signature $\sigma$ with support $\cS$ is said to be associated with a function $g \in C(\Omega)$ if $\cS$ is included in the set  $\{\omega \in \Omega: |g(\omega)| = \|g\|_\Omega \}$ of extremal points of~$g$  and if $\sigma(\omega) = \sgn(g(\omega))$ for all $\omega \in \cS$.
\edefn
 
The argument outlined before the definition \rev{justifies} the following brief statement, found e.g. in \cite[Theorem 2.6]{Riv}. 
We emphasize that it does not provide a way to find extremal signatures and best approximants,
but if one comes up with a guess for these (as we will do in Section \ref{SecEx}),
then it provides a way to verify that the guess is correct.

\bthm
\label{ThmSig}
An element $v^* \in \cV$ is a best approximant to $f \in C(\Omega)$ from $\cV$
if and only if there exists an extremal signature $\sigma$ for $\cV$ associated with $f-v^*$.
Moreover, the support of such a signature can be chosen to have size $\le \dim(\cV)+1$.
\ethm

An important detail not made apparent in the above statement is the existence of a signature common to all best approximants---this is revealed by \eqref{IntroSigna},
because the involved arguments did not depend on the best approximant $v^*$.
This fact explains our preference for solving \eqref{DualForm} over~\eqref{PrimalForm},
especially since \eqref{PrimalForm} typically have nonunique solutions.

\subsection{Simple reductions}

Given a fixed number $d$ of variables and a fixed degree $n$,
completely solving the problem of $d$-variate $n$th degree Chebyshev polynomials
requires finding best approximants to all $d$-variate $n$th degree monomials,
so we would a priori need to tackle $\binom{n+d-1}{d-1}$ subproblems.
For $d=3$ and $n=6$, it amounts to $28$ subproblems
and for $d=3$ and $n=10$, it amounts to $66$ subproblems.
Fortunately, this number can be decreased drastically by leveraging two simple reductions.
The first reduction allows us to discard the indices $k_i=0$ in the multiindex $k$ of the monomial $m_k$,
provided the full problem has been solved for all $d'<d$.
The second reduction allows us to consider only indices $k_1,\ldots,k_d$ that are ordered from largest to smallest, say.
These facts are precisely stated in the two propositions below.
In the first one,
the domain $\Omega$ can be taken as any of our preferred choices---the simplex $S$,
the cross-polytope $C$,
the euclidean ball $B$,
or the hypercube $H$---by selecting $\varphi = 0$.
The argument, already found in \cite[Proposition 4.1]{X05} for $\Omega = B$, is included here to also cover the case $\varphi \not= 0$.
\rev{The statement uses the notation
$\bN_0$ for $\{0,1,2,\ldots\}$,
$\omega_I \in \bR^I$ for the restriction of a vector $\omega \in \bR^d$ to a subset $I$ of $\{1,2,\ldots,d\}$,
and $I^c$ for the complement of $I$.}

\bprop
\label{PropRedu1}
Given $k \in \bN_0^d$ with $|k|=k_1+\cdots+k_d=n$,
let $I \coloneqq \left\{ i = 1,\ldots,d: k_i > 0 \right\}$,
\revv{let $d' \coloneqq |I|$, and let $k' \coloneqq k_I \in \bN_0^{d'}$, 
which satisfies $|k'|=k'_1+\cdots+k'_{d'}=n$.}
Let \revv{also} $\Omega' \inc \bR^I$ be the $d'$-dimensional domain defined by 
$\Omega' = \{ \omega_I, \, \omega \in \Omega \}$.
\revv{Suppose} that there is a $\varphi \in \bR^{I^c}$ such that the element $\wt{\omega}$ defined by $\wt{\omega}_i = \omega_i$ for $i \in I$
and $\wt{\omega}_i = \varphi_i$ for $i \in I^c$
belongs to $\Omega$ whenever $\omega \in \Omega$.
\revv{Then one has}
$$
E_{\cP_{n-1}^d}(m_k,\Omega) 
= E_{\cP_{n-1}^{d'}}(m_{k'},\Omega').
$$
\eprop

\bpf
On the one hand, with
$q' \in \cP_{n-1}^{d'}$ such that
$E_{\cP_{n-1}^{d'}}(m_{k'},\Omega') = \|m_{k'} - q' \|_{\Omega'}$
and with  $q \in \cP_{n-1}^d$ defined by $q(x) = q'(x_I)$, $x \in \bR^d$,
we have
\begin{align*}
E_{\cP_{n-1}^{d'}}(m_{k'},\Omega')
& = \|m_{k'} - q' \|_{\Omega'}
= \max_{\omega \in \Omega} |m_{k'}(\omega_I) - q'(\omega_I)|
= \max_{\omega \in \Omega} |m_{k}(\omega) - q(\omega)|
= \|m_k - q\|_\Omega\\
& \ge 
E_{\cP_{n-1}^d}(m_k,\Omega).
\end{align*}   
This inequality was obtained independently of the existence of $\varphi$.
On the other hand,
for the reverse inequality,
given $p \in \cP_{n-1}^d$,
we have 
$$
\|m_k - p\|_\Omega
= \max_{\omega \in \Omega} |m_k(\omega) - p(\omega)|
\ge \max_{\omega \in \Omega} |m_{k}(\wt{\omega}) - p(\wt{\omega}) |
= \max_{\omega \in \Omega} |m_{k'}(\omega_I) - \wt{p}(\omega_I) |,
$$
where $\wt{p}$ is implicitly defined as a polynomial in $\cP_{n-1}^{d'}$.  
Therefore, we obtain
$$
\|m_k - p\|_\Omega \ge \|m_{k'} - \wt{p}\|_{\Omega'}
\ge E_{\cP_{n-1}^{d'}}(m_{k'},\Omega').
$$
The inequality $E_{\cP_{n-1}^d}(m_k,\Omega) \ge E_{\cP_{n-1}^{d'}}(m_{k'},\Omega')$ follows by taking the infimum over $p$.
\epf

The second fact, stated hereafter,
has been previously used to derive a number of examples in \cite{AV1, AV2, X04}.
We include a standard argument (see e.g. \cite[Theorem 3.2]{X05}) for the convenience of the reader.
This fact is to be used with $\cV = \cP_{n-1}^d$ and $G$ being the group of permutation\revv{s} of~$\{1,2,\ldots,d\}$.

\bprop
\label{PropRedu2}
Given a finite group $G$ acting on a domain $\Omega \inc \bR^d$,
for $h \in C(\Omega)$ and $g \in G$,
let  $h_g \in C(\Omega)$ be defined by $h_g(\omega) = h(\rev{g \omega})$, $\omega \in \Omega$.
If the domain $\Omega$ and the subspace $\cV \inc C(\Omega)$ are invariant under the action of $G$,
in the sense that
\begin{align*}
\Omega_g & \coloneqq \{\rev{g \omega }, \omega \in \Omega \} & & \mbox{coincides with $\Omega$ for all } g \in G,\\
\cV_g & \coloneqq \{ v_g, v \in \cV \} & & \mbox{coincides with $\cV$ for all } g \in G,
\end{align*}
then, for any $f \in C(\Omega)$ and any $g \in G$,
$$
E_\cV(f,\Omega) = E_\cV(f_g,\Omega).
$$
Furthermore, if $f$ is invariant under the action of $G$, i.e., if $f_g$ coincides with $f$ for all $g \in G$, then there is a best approximant $v^*$ to $f$ from $\cV$ which is  invariant under the action of $G$,
i.e., $v^*_g = v^*$ for all $g \in G$. 
\eprop

\bpf
For $f \in C(\Omega)$,
let $v' \in \cV$ be a best approximant to $f$ from $\cV$.
The invariance of $\Omega$ implies that,
for any $g \in G$,
\begin{align*}
E_\cV(f,\Omega) &  = \max_{\omega \in \Omega} |f(\omega)-v'(\omega)| = \max_{\omega \in \Omega} |f(\rev{g \omega})-v'(\rev{g \omega})|
= \max_{\omega \in \Omega} |f_g(\omega)-v'_g(\omega)| \\
& \ge E_\cV(f_g,\Omega),
\end{align*}
where the last step relied on the invariance of $\cV$ to ensure that $v'_g \in \cV$.
A similar argument with $f_g$ in place of $f$ and $g^{-1}$ in place of $g$ would yield the reverse inequality $E_\cV(f_g,\Omega) \ge E_\cV(f,\Omega)$ and in turn the desired equality.
Now, let us assume in addition that $f_g = f$ for all $g \in G$.
Using the above,
we have $\|f - v'\|_\Omega = \|f - v'_g \|_\Omega$,
so that $v'_g$ is also a best approximant to $f$ from $\cV$ for all $g \in G$.
Consequently, the element $v^* \coloneqq |G|^{-1} \sum_{g \in G} v'_g \in \cV$,
as a convex combination of best approximants,
is itself a best approximant to $f$ from $\cV$.
It is also readily seen that $v^*$ thus defined is invariant under the action of $G$.
\epf

As alluded to before,
Propositions \ref{PropRedu1} and \ref{PropRedu2} imply that,
for the simplex, the cross-polytope, the euclidean ball,
and the hypercube, 
it is enough to consider the monomials $m_k$
where $k \in \bN_0^d$ satisfies $k_1+\cdots+k_d=n$ and $k_1 \ge \cdots \ge k_d \ge 1$.
The number of these monomials equals the number of partitions of the integer $n$ into exactly $d$ parts.
This number $p_d(n)$ is known to obey the recurrence relation $p_d(n) = p_{d-1}(n-1) + p_d(n-d)$,
which allows one to arrange them in a triangular table akin to Pascal's triangle.
For instance, for  $d=3$ and $n=10$, one has $p_3(10)=8$.
For $d=3$ and $n=6$, one has $p_3(6)=3$, with the three partitions being $(4,1,1)$, $(3,2,1)$, and $(2,2,2)$.
In case of the cross-polytope,
the values of the three corresponding errors of best approximation are reported in the last column of Table \ref{TableCrossPolytope}, none of which were known before.

\subsection{Known multivariate Chebyshev polynomials}

In this section, we gather previously obtained results about multivariate Chebyshev polynomials for our domains of interest, with the exception of the cross-polytope,
which seems to have been cast aside in the literature.
We will use from now on the shorthand notation
$$
E(k,\Omega) \coloneqq E_{\cP_{n-1}^d}(m_k,\Omega),
$$ 
since considering $k=(k_1,\ldots,k_d)$ with $k_1+\cdots+k_d=n$ implicitly tells us the value of $d$ and $n$.

\paragraph{The hypercube.} 

The case of the hypercube, i.e., $\Omega = H$, is completely resolved. 
Indeed, the geometry of the domain bodes well for calculations involving the tensor products of univariate polynomials $p_1,\ldots,p_d$, as defined by $(p_1 \otimes \cdots \otimes p_d)(x_1,\ldots, x_d) = p_1(x_1) \cdots p_d(x_d)$.
It is not difficult to establish the following  result by invoking  signatures.

\bthm
Given $k \in \bN^d$ with $k_1+\cdots+k_d=n$,  one has
$$
E(k,H) = 2^{-n+d}
$$
and a best approximant to $m_k$ from $\cP_{n-1}^d$ is given by $m_k - 2^{-n+d} T_{k_1} \otimes \cdots \otimes T_{k_d}$.
\ethm

This result was proved by several authors,
see e.g. \cite{Slo,EhlZel}.
In \cite{VY}, it has also been shown that $m_k - 2^{-n+d} T_{k_1} \otimes \cdots \otimes T_{k_d}$ is a unique best approximant when and only when $d = 2$ and $k_1 = k_2$.

\paragraph{The euclidean ball.} 

The case of the euclidean ball, i.e., $\Omega = B$, is partially resolved: it is solved for $d=2$ variables but not completely in $d \ge 3$ variables.
With $\rev{U_\ell  = T_{\ell+1}'/(\ell+1)}$ denoting the univariate $\ell$th degree Chebyshev polynomial of the second kind,
the result for $d=2$ reads as follows.

\bthm
Given $k \in \bN^2$ with $k_1+k_2=n$, one has
$$
E(k,B) = 2^{-n+1}
$$
and a best approximant to $m_k$ from $\cP_{n-1}^2$ is given by $m_k - 2^{-n} (U_{k_1} \otimes U_{k_2} - U_{k_2-2}\otimes U_{k_1-2})$,
with the understanding that $U_{-1}=0$.
\ethm

This was obtained for the first time in  \cite{G}. 
Other explicit best approximants can be found in~\cite{BHN, R}. 
It is also known that the difference between two best approximants 
has the form $(1-x_1^2-x_2^2) q(x_1,x_2)$ for some $q \in \cP_{n-3}^2$.

For $d > 2$,  best approximants to monomials are known only for a few low-degree instances,
such as $m_{(1,\ldots,1)}(x) = x_1 x_2\cdots x_d$ and $m_{(2,1,\ldots,1)}(x) = x_1^2 x_2\cdots x_d$, see \cite{AV1, AV2, X04}. 
Restricting our attention to the case $d =3$,  we now cite some articles and the result they contain:\footnote{The value of $E((3,1,1),B)$ is implicit in \cite[Theorem 2.1.(b)]{MoaPeh}---deriving the explicit value requires some work.}$^{,}$\footnote{There was a typographical error concerning the value of $E((4,4,4),B)$
in \cite[Theorem 3.2]{X04}.}
\begin{align*}
\cite{AV1} : & & E((1,1,1),B) & = 3^{-3/2},  \\
\cite{AV2} : & & E((2,1,1),B) & = (3-\sqrt{8})/2,\\
\cite{MoaPeh} : & &  E( (3,1,1),B) & = (1-a)(a^3/5)^{1/4}/5,
\quad a= \mbox{smallest root of } 9t^4 - 29 t^3 + 24 t^2 -29 t +9,\\ 
 \cite{X04} : & &
E((2,2,2),B) & = 1/72,\\
\cite{X04}: & &
E((4,4,4),B) & = b^{-1}/ 27^2,  \qquad \qquad \quad \;\,\, b \approx 21.8935834.
\end{align*}

\paragraph{The simplex.} 

The case of the simplex, i.e., $\Omega = S$, is also partially resolved.
Indeed, for $d=2$,  best approximants to monomials are presented in \cite{NewXu}.
The result is recalled below.

\bthm
Given $k \in \bN^2$ with $k_1+k_2=n$,
one has
$$
E(k,S) = 2^{-2n+1}
$$
and a best approximant to $m_k$ from $\cP_{n-1}^2$ is given by $m_k-T_{k_1,k_2}$,
where
$$
T_{k_1,k_2}(x,y) = T_{k_1-k_2}(2x-1)T_{k_2}(8xy-1)+8xy(2x-1)U_{k_1-k_2-1}(2x-1)U_{k_2-1}(8xy-1)
$$
for $k_1 \ge k_2$, with the understanding that $U_{-1}=0$.
\ethm
In the case $d=3$, we mention the results
\begin{align*}
\cite{X04}: & &
E((1,1,1),S) & = 1/72,\\
\cite{X04}: & &
E((2,2,2),S) & = b^{-1}/ 27^2,  \qquad \quad b \approx 21.8935834.
\end{align*}
Note that there is a close connection between the best approximants  on the simplex $S$ to the monomial $m_k(x)=x_1^{k_1} \cdots x_k^{k_d}$  and the best approximants on the euclidean ball $B$ to the monomial $m_{2k}(x)=x_1^{2k_1} \cdots x_k^{2k_d}$,
see \cite{X04} for the precise statement.

\section{Description of the Computational Procedure}
\label{SecDes}

In this section,
we explain the procedure that we derived and exploited in order to produce a number of new multivariate Chebyshev polynomials (uncovered in Section \ref{SecEx}).
\revv{Our implementation is
available at \url{https://github.com/foucart/Multivariate_Chebyshev_Polynomials}.}
\revv{Although it} limits itself to the best approximation from $\cP_{n-1}^d$ to monomials $m_k \in \cP_n^d$ on the hypercube, the euclidean ball, the cross-polytope, and the simplex,
the underpinning procedure is more general.
\revv{Indeed, }
it could handle any polynomial $f \in \cP_{N}^d$, $N \ge n$, instead of $m_k$,
while the domain $\Omega \in \bR^d$ could be any \rev{compact} (with nonempty interior)
basic semialgebraic set.
\revv{This means}
that there exist polynomials $g_1,\ldots,g_H$ such that
\be
\label{Om}
\Omega = \{ x \in \bR^d:  g_1(x) \ge 0, \ldots,g_H(x) \ge 0 \}.
\ee
\revv{We further assume that the polynomials describing $\Omega$ satisfy} the Archimedean condition,
meaning e.g. that 
there exist a constant $C>0$ and sum-of-squares polynomials $\sigma_0,\sigma_1,\ldots,\sigma_H$ such that
$$
C - \|x\|_2^2 = \sigma_0(x) + \sum\nolimits_{h=1}^H \revv{ \sigma_h(x) g_h(x)}
\qquad 
\mbox{for all } x \in \Omega.
$$
\revv{This occurs for all four domains considered in this article}---for instance,
the argument for the simplex $S$ can be found in \cite[Example 12.49]{Schm} and a similar argument applies to any compact convex polytope.

\rev{Our strategy to deal with the best approximation problem is to transform it into} an instance of the Generalized Moment Problem (GMP), so that we can invoke the Moment-SOS hierarchy designed to solve a GMP whose data are algebraic 
(polynomials and semialgebraic sets), see \cite{JBL}. 
This process leverages a combination of:
(i)~\emph{semidefinite programming}\footnote{A semidefinite program (SDP) is a conic convex optimization problem which can be solved in time polynomial in its input size, up to arbitrary (fixed) precision (e.g. with interior point methods); for more details, see e.g. \cite{LauRen}.},
an efficient machinery in Convex Optimization developed since the late seventies, 
and 
(ii)~powerful positivity certificates and their dual analogs concerning the moment problem.
These two ingredients were not available at
the time of pioneering works such as the paper \cite{shapiro} by Rivlin and Shapiro,
in which dual formulations were mostly used to prove the existence of optimal solutions and to characterize them.
For the numerical computations, we will rely on {\sf GloptiPoly 3},
since many of the GMP components are built in this {\sc matlab}/{\sc octave} program, see \cite{GloptiPoly}.

Let us be more specific about the computation of the error of best approximation by polynomials from $\cP_{n-1}^d$ to a polynomial $f \in \cP_{N}^d$, $N \ge n$,
i.e., of $E^*\coloneqq E_{\cP_{n-1}^d}(f,\Omega)$,
\revv{viewed as the optimal value of \eqref{PrimalForm} or of \eqref{DualForm}---which we call primal program and dual program, respectively. 
Looking at \eqref{PrimalForm} first, it can be reformulated along the lines of \eqref{ReformNonNegCst} into
$$
E^* = \min_{c \in \bR, p \in \cP_{n-1}^d} \; c
\qquad \mbox{subject to }  \left\{\begin{matrix}  c + f- p \in \cP_{+}^d(\Omega),\\
c - f + p \in \cP_{+}^d(\Omega),  \end{matrix} \right.
$$
where $\cP_{+}^d(\Omega)$ denotes the cone of $d$-variate polynomials that are nonnegative on $\Omega$.
As for \eqref{DualForm}, 
\revv{by identifying linear functionals $\la \in C(\Omega)^*$  with finite signed Borel measures $\mu$ on $\Omega$ through 
$\la(g) = \int_\Omega g d\mu$, $g \in C(\Omega)$,
it can be rewritten as
\be
\label{InfDual0}
E^* = \max_{\mu \in \mathcal{M}(\Omega)} 
\int_\Omega f \, d\mu
\quad \mbox{s.to } \quad \int_\Omega p \, d\mu = 0 \mbox{ for all }
p \in \cP_{n-1}^d \quad \mbox{and } \int_\Omega d|\mu|  = 1.
\ee
With $\mathcal{M}_+(\Omega)$ denoting the cone of finite nonnegative Borel measures on $\Omega$,
the latter can further be rewritten as\footnote{\revv{\eqref{InfDual0}$\le$\eqref{InfDual} holds by considering the Jordan decomposition $\mu =: \mu^+-\mu^-$ of a $\mu \in \cM(\Omega)$ which is optimal for \eqref{InfDual0}
and \eqref{InfDual}$\le$\eqref{InfDual0} holds by considering $\mu := (\mu^+-\mu^-) / \int_\Omega d|\mu^+-\mu^-|$ 
where $\mu^+,\mu^- \in \cM_+(\Omega)$
are optimal for \eqref{InfDual}.}}
\begin{align}
\label{InfDual}
E^* = \max_{\mu^\pm \in \mathcal{M}_+(\Omega)} 
\int_\Omega f \, d(\mu^+-\mu^-)
&\quad \mbox{s.to } \int_\Omega p \, d(\mu^+-\mu^-) = 0 \mbox{ for all }
p \in \cP_{n-1}^d\\
\nonumber
& \quad \mbox{and } \int_\Omega d(\mu^+ + \mu^-) = 1.
\end{align}}This is an instance of the  Generalized Moment Problem (GMP), mentioned above.
Note that $\cP_{+}^d(\Omega)$ is the dual cone of $\mathcal{M}_+(\Omega)$ (and vice versa) with respect to the pairing
$
\langle p,\mu\rangle := \int_\Omega p d\mu.
$
In other words, one has $\mathcal{M}_+(\Omega)^* = \cP_{+}^d(\Omega)$, where the dual cone is defined in the usual way, i.e.,
\[
\mathcal{M}_+(\Omega)^* := \left\{ p \in \cP^d(\Omega) \; : \; \langle p, \nu\rangle \ge 0  \mbox{ for all } \nu \in \mathcal{M}_+(\Omega)\right\}.
\]
}\revv{The {\em quadratic module} of  $g_1,\ldots,g_H$,
denoted by $\mathcal{Q}(g_1,\ldots,g_H)$,
is an important subcone of $\cP_{+}^d(\Omega)$.
It consists of all polynomials of the form
\begin{equation}
\label{eq:quadratic module}
p\,=\,\sigma_0+\sigma_1\,g_1+\cdots +\sigma_H\,g_H\,
\end{equation}
for some sum-of-squares polynomials $\sigma_0,\sigma_1,\ldots,\sigma_H$.
The associated {\em truncated } quadratic module of degree $2t$, denoted by $\mathcal{Q}(g_1,\ldots,g_H)_t$, 
is the subset of $\mathcal{Q}(g_1,\ldots,g_H)$ 
obtained by restricting all the terms in \eqref{eq:quadratic module} to have degree at most $2t$, 
i.e., $\deg(\sigma_0) \le 2t$ and $\deg(\sigma_h g_h) \le 2t$ for $h = 1,\ldots,H$.
}

\revv{
By definition, the dual cone of $\mathcal{Q}(g_1,\ldots,g_H)$ consists of linear functionals $\la$ on $\cP^d$  
such that $\la(p) \ge 0$ for all $p \in \mathcal{Q}(g_1,\ldots,g_H)$.
This dual cone can be described in terms of sequences of positive semidefinite  matrices.
\revv{
To this end, we denote by ${\rm Hank}_t(y)$  the Hankel matrices of size $ \binom{t+d}{d} \times \binom{t+d}{d}$ built from a sequence $y$ indexed by~$\bN_0^d$
as having entries $y_{i + j}$ for $i,j \in \{0,\ldots,t\}^d$.
 Similarly,
we construct matrices ${\rm Hank}_t(G_h y)$
from sequence $G_h y$ associated with each $g_h$,
namely $(G_h y)_\ell = \sum\limits_{|\ell'| \le \deg(g_h)} {\rm coeff}_{m_{\ell'}}(g_h) y_{\ell + \ell'}$ for any $\ell \in \bN_0^d$.}
}

\revv{
The following theorem by Putinar \cite{Put} describes the relation between $\cP_{+}^d(\Omega)$ and $\mathcal{Q}(g_1,\ldots,g_H)$, 
and, on the dual side, between $\mathcal{M}_+(\Omega)$ and $\mathcal{Q}(g_1,\ldots,g_H)^*$.

\bthm
\label{th-put}
Let $\Omega = \{x \in \bR^d: g_1(x) \ge 0,\ldots, g_H(x) \ge 0 \}$ be a compact semialgebraic set 
such that $g_1,\ldots,g_H$ satisfy an Archimedean condition. Then:\vspace{2mm}\\
(Real analysis facet):
A real sequence $y=(y_\ell)_{\ell\in\N^d_0}$ has a representing measure on $\Omega$ if and only if
$$
{\rm Hank}_t(y)\succeq 0, {\rm Hank}_t(G_1\,y)\succeq 0,\ldots,{\rm Hank}_t(G_H\,y)\succeq 0
\quad \mbox{for all } t\in\N .
$$
Thus, a linear functional $\la$ on $\cP^d$ is (identified with an element) in $\mathcal{M}_+(\Omega)$ if and only if
 the sequence $ y = \left( \la(m_\ell)\right)_{\ell \in \bN_0^d}$ satisfies the above  semidefinite conditions.\vspace{2mm}\\
(Real algebraic geometry facet):  If  $p\in \cP^d$ is strictly positive on $\Omega$, then $p \in \mathcal{Q}(g_1,\ldots,g_H)$.
Conversely, any $p \in \mathcal{Q}(g_1,\ldots,g_H)$ is nonnegative on $\Omega$.
\ethm
}

\revv{Putinar's Positivstellensatz may also be formulated as follows:
if a polynomial $p \in \cP^d$ is strictly positive on $\Omega$, 
then there exists $t \in \mathbb{N}$ such that $p \in  \mathcal{Q}(g_1,\ldots,g_H)_t$.
As a result, one has
\begin{align}
\label{added_explanation}
  E^* 
 &= \inf_{c \in \bR, p \in \cP_{n-1}^d} \left\{c \, : \, c \pm (f-p)(x) >0 \mbox{ for all } x\in \Omega  \right\} \\
 \nonumber 
 &= \inf_{c \in \bR, p \in \cP_{n-1}^d} \left\{c \, : \, c \pm (f-p) \in  \mathcal{Q}(g_1,\ldots,g_H)_t \mbox{ for some } t \in \bN \right\}.
\end{align}
Since sums-of-squares can be modeled using semidefinite matrices,
semidefinite programming can be used to compute,
for any fixed $t \in \bN$ large enough\footnote{\revv{The semidefinite program of \eqref{DefLB} may be infeasible for small $t$'s, in which case on has ${\rm ub}_t = \infty$,
but Putinar's Positivstellensatz guarantees that it is feasible for sufficiently large $t$.}},
an upper bound for $E^*$ as
\begin{equation}
\label{DefLB}
{\rm ub}_t := \inf_{c \in \bR, p \in \cP_{n-1}^d} \; c
\qquad \mbox{subject to }  \left\{\begin{matrix}  c + f- p \in \mathcal{Q}(g_1,\ldots,g_H)_t,\\
c - f + p \in \mathcal{Q}(g_1,\ldots,g_H)_t.  \end{matrix} \right.
\end{equation}
Notice that $({\rm ub}_t)$ forms a nonincreasing sequence and,
according to \eqref{added_explanation},
that $\lim_{t \to \infty} {\rm ub}_t = E^*$.
The idea to construct such a sequence of semidefinite programs is originally due to Lasserre \cite{lass-siopt-01}.
}

\revv{The sequence of dual semidefinite programs leads us to consider,
for each  $t \in \mathbb{N}$, 
\begin{align}
\label{InfDual_t}
{\rm ub}'_t
:= \sup_{y^{\pm} \in \bR^{\{0,\ldots,2t\}^d}} & \;
\sum_{|\ell| \le N} {\rm coeff}_{m_\ell}(f) (y^+_\ell - y^-_\ell)
 \quad \mbox{s.to } y^+_\ell - y^-_\ell = 0, \, |\ell| \le n-1, \; \, y^+_0 + y^-_0 = 1, \\ 
 \nonumber
 \mbox{and } \; \; \; \; & \quad {\rm Hank}_t(y^{\pm}) \succeq 0,
{\rm Hank}_{t-n'_1}(G_1 y^{\pm}) \succeq 0,
\ldots, 
{\rm Hank}_{t-n'_H}(G_H y^{\pm}) \succeq 0,
\end{align}
where $n'_1 := \lceil \mbox{deg}(g_1)/2\rceil, \ldots, n'_H := \lceil \mbox{deg}(g_H)/2\rceil$.
This program can be thought of as a truncation at level $t$ of the program \eqref{InfDual} in which $\mu^\pm \in \cM_+(\Omega)$ would be replaced by their infinite sequences of moments $y^\pm_\ell = \int_\Omega m_\ell \, d\mu^\pm$, $\ell \in \bN_0^d$.
Here, the finite sequence $y^\pm$ represents pseudo-moments.
Since the feasibility set of \eqref{InfDual} is contained in the feasibility set of \eqref{InfDual_t},
one has $E^* \le {\rm ub}'_t$,
i.e., the ${\rm ub}'_t$ are also upper bounds for $E^*$.
Furthermore,
by weak duality for semidefinite programs, one always has ${\rm ub}'_t \le {\rm ub}_t$.
In view of $\lim_{t \to \infty} {\rm ub}_t =~E^*$,
it is therefore immediate  that 
$\lim_{t \to \infty} {\rm ub}'_t = E^*$, too.
The strong duality result ${\rm ub}'_t = {\rm ub}_t$ 
is actually valid,
since a constraint qualification always holds, by a similar argument to \cite[Theorem 4.2(a)]{lass-siopt-01}.
For instance, one may construct a strictly feasible solution to the dual semidefinite program
by setting $\mu^+ = \mu^-$ equal to the uniform measure on $\Omega$. 
Having assumed that $\Omega$ has a nonempty interior, 
this implies that $\mu^+$ and $\mu^-$ lie in the interior of $(\mathcal{Q}(g_1,\ldots,g_H)_t)^*$. 
The interested 
reader can find more details in the proof of \cite[Theorem 4.2(a)]{lass-siopt-01}.
In particular, the dual semidefinite program is unbounded if the primal program \eqref{DefLB} is infeasible.

In summary, the following result constituting the basis for our computational procedure has been established.

\bthm
\label{ThmHiera}
Let $\Omega = \{x \in \bR^d: g_1(x) \ge 0,\ldots, g_H(x) \ge 0 \}$ be a compact semialgebraic set 
such that $g_1,\ldots,g_H$ satisfy an Archimedean condition
and let $f \in \cP_{N}^d$  a polynomial $f$ of degree $N \ge n$.
The values ${\rm ub}_t$ and ${\rm ub}'_t$ computed in~\eqref{DefLB} and \eqref{InfDual_t}
with $t \ge N + \max\{n'_1,\ldots,n'_H\}$, $n'_h\coloneqq \lceil \deg(g_h)/2 \rceil$,
form nonincreasing sequences of upper bounds for the error $E_{\cP_{n-1}^d}(f,\Omega)$ of best approximation 
to~$f$ from $\cP_{n-1}^d$ relative to $\Omega$.
Moreover, these sequences converge to the error of best approximation,
i.e.,
$$
\lim_{t \to \infty} {\rm ub}_t 
= \lim_{t \to \infty} {\rm ub}'_t
= E_{\cP_{n-1}^d}(f,\Omega).
$$
\ethm

}

\revv{The above ideas have been implemented in the commands {\sf ChebPoly\_primal} and {\sf ChebPoly\_dual},
which require {\sf GloptiPoly~3}.
Although any semidefinite solver could have been used to solve \eqref{DefLB} and \eqref{InfDual_t}, 
we chose {\sf GloptiPoly~3}
because it offers a more user-friendly way to formulate the problems.
Another advantage is that,
when solving \eqref{InfDual_t},
{\sf GloptiPoly~3} not only outputs moments
but it can also extract the atomic measure generating them---the procedure used to do so is described in \cite{GloptiPoly solution extraction}.
This corresponds to situations where the convergence of the sequence of upper bounds occurs in a finite number of steps, which is often the case in practice. 
While we do not know this number of steps in advance, 
in such cases, the precise value of $E_{\cP_{n-1}^d}(f,\Omega)$
can be obtained by computing ${\rm ub}'_t$ (for $t$ large enough)
and certifying its genuine optimality ``numerically'' within {\sf GloptiPoly~3}.
Indeed, based on a deep result of Curto and Fialkow, 
{\sf GloptiPoly~3} verifies a condition ensuring that 
${\rm ub}'_t = E_{\cP_{n-1}^d}(f,\Omega)$ through a numerical estimation a rank drop in the Hankel matrices,
see~\cite{GloptiPoly solution extraction}. 
Such certificates 
substantiate the numerical values presented in the next section.}

\section{New Explicit Results}
\label{SecEx}

In this section,
we present the discoveries that were made by exploiting the computational procedure just described.
In particular, we consider the situation of $d=3$ variables and give the values of the errors of best approximation from $\cP_{n-1}^3$ to (essentially) all monomials $m_k$ of degree $n$ up to~$6$.
Excluding the already settled case where $\Omega$ is the hypercube,
we report,
from largest to smallest~$\Omega$,
on the cases of the euclidean ball, the cross-polytope, and the simplex in Tables \ref{TableBall}, \ref{TableCrossPolytope}, and \ref{TableSimplex}, respectively.
\revv{Their content can be reproduced by running the codes made available on \url{https://github.com/foucart/Multivariate_Chebyshev_Polynomials}. Note that an asterisk ($\ast$) appended to a value means that it could not be numerically certified\footnote{\revv{To be sure, it is the value of the error of best approximation that is certified,
not the extracted atomic measure.
In the case of $m_{(2,2,1)}$ on the simplex, for instance, even though the value $E((2,2,1),S) \approx 4.695 \times 10^{-4}$ is certified, the extracted atomic measure suffers from numerical inaccuracies. 
This case is challenging---we have tried to derive an explicit solutions from our computations, without any success.}} by {\sf GloptiPoly~3}.}
In a few instances, we could distill explicit expressions for previously unknown Chebyshev polynomials, \rev{see Theorems \ref{ThmB221} and \ref{ThmS211}}.

\subsection{The euclidean ball}

When $\Omega = B$,
the values of $E(k,\Omega)$ had earlier been found when $|k|=3$ and $|k|=4$,
but not for $|k|=5$ (except $k=(3,1,1)$) nor $|k|=6$ (except $k=(2,2,2)$).
All these values are shown in Table~\ref{TableBall}.
Here and in other tables,
the value of $E(k,\Omega)$ needs to be multiplied by the factor presented at the top of its column,
so e.g. $E((4,1,1),B) \approx 1.923 \times 10^{-2}$.

\begin{table}[h]
\small
\begin{center}
\begin{tabular}{|c|c|c|c|}
\hline
degree $n=3$ ($\times 10^{-1}$) & degree $n=4$ ($\times 10^{-2}$) & degree $n=5$ ($\times 10^{-2}$) & degree $n=6$ ($\times 10^{-2}$)\\
\hline
\hline
$E((1,1,1),B) \approx 1.924$ & $E((2,1,1),B) \approx 8.578$ & $E((3,1,1),B) \approx 4.016$ & $\mathbf{E((4,1,1),B) \approx 1.923}$ \\
\hline
 & & $\mathbf{E((2,2,1),B) \approx 3.630}$ & $\mathbf{E((3,2,1),B) \approx 1.652}$\\
\hline
 & & & $E((2,2,2),B) \approx 1.388$\\
 \hline 
\end{tabular}
\end{center}
\caption{Euclidean ball in dimension $d=3$: the previously unknown values are shown in boldface.}
\label{TableBall}
\end{table}

In the case $k=(2,2,1)$,
it was possible to recognize (part of) 
the signature points,
which led us to deriving a Chebyshev polynomial explicitly.
The result reads as follows.

\bthm
\label{ThmB221}
With $a \coloneqq\max \big\{ (1+t)^2 (1-t) t/(4(1+4t+4t^2)) , t \in [0,1] \big\} \approx 3.63000825 \times 10^{-2}$,
the error of best approximation on the euclidean ball to $m_{(2,2,1)}(x_1,x_2,x_2) = x_1^2 x_2^2 x_3$ by trivariate polynomials of degree at most $4$ is
$$
E_{\cP_4^3}(m_{(2,2,1)},B) = a,
$$
while a Chebyshev polynomial is given by 
$$
P(x) = m_{(2,2,1)}(x_1,x_2,\revv{x_3}) + a \, T_3(x_3),
$$
where $T_3(t) = 4 t^3 - 3t$ is the univariate Chebyshev polynomial of degree $3$.
\ethm

\bpf
Guided by our computational procedure,
we make the guess---which we are about to verify---that a signature has support $\cS = \cS^+  \cup \cS^-$,
where $\cS^- = - \cS^+$ and 
\begin{multline*}
\cS^+ = \left\{
\bpmx 0 \\ 0 \\ 1 \epmx,
\bpmx \sqrt{3}/2 \\ 0 \\ -1/2 \epmx,
\bpmx -\sqrt{3}/2 \\ 0 \\ -1/2 \epmx,
\bpmx 0 \\ \sqrt{3}/2 \\ -1/2 \epmx,
\bpmx 0 \\ -\sqrt{3}/2 \\ -1/2 \epmx,\right.
\\
\left.
\bpmx \sqrt{(1-\tau^2)/2} \\ \sqrt{(1-\tau^2)/2} \\ \tau \epmx,
\bpmx \sqrt{(1-\tau^2)/2} \\ -\sqrt{(1-\tau^2)/2} \\ \tau \epmx,
\bpmx -\sqrt{(1-\tau^2)/2} \\ -\sqrt{(1-\tau^2)/2} \\ \tau \epmx,
\bpmx \sqrt{(1-\tau^2)/2} \\ -\sqrt{(1-\tau^2)/2} \\ \tau \epmx
\right\}.
\end{multline*}
For reasons soon to become apparent,
the parameter $\tau \approx 0.4052$ is chosen as the maximizer of $(1+t)^2 (1-t) t /(4(1+4t+4t^2))$ over $t \in [0,1]$,
so that $(1+\tau)^2 (1-\tau) \tau/4 = a (1+4 \tau+ 4 \tau^2)$.
Multiplying throughout by $(1-\tau)$ yields
$((1-\tau^2)/2)^2 \tau = a (1-T_3(\tau))$.
As for best approximants to $m_{(2,2,1)}$ from $\cP_4^3$,
we shall show,
in two steps, that $p^*(x_1,x_2,x_3) \coloneqq -a T_3(x_3)$ is one of them.

The first step consists in proving
that $| (m_{(2,2,1)} - p^*)(x)| = \| m_{(2,2,1)} - p^* \|_B$ for all $x \in \cS$.
To see this, we start by noticing that $(m_{(2,2,1)} - p^*)(x) = a$ for all $x \in \cS^+$---this is a simple verification by plugging in the values $x \in \cS^+$ into $(m_{(2,2,1)} - p^*)(x)$,
but we emphasize that $(m_{(2,2,1)} - p^*)(\pm \sqrt{1-\tau^2},\pm \sqrt{1-\tau^2}, \tau) = ((1-\tau^2)/2)^2 \tau + a T_3(\tau) = a$ owes to our choice of $\tau$.
Then,
from $\cS^- = - \cS^+$ and the oddity of $m_{(2,2,1)}$,
it follows that $(m_{(2,2,1)} - p^*)(x) = -a$ for all $x \in \cS^-$.
All in all, we arrived at $|(m_{(2,2,1)} - p^*)(x)| = a$ for all $x \in \cS$.
Next, we claim that $|(m_{(2,2,1)} - p^*)(x)| \le a$ for all $x \in B$.
By the oddity of $m_{(2,2,1)}$ again,
it is enough to establish, for all $x \in B$, that
$(m_{(2,2,1)} - p^*)(x) \le a$,
i.e., that $x_1^2 x_2^2 x_3 \le a (1-T_3(x_3))$.
If $x_3 \in [-1,0]$, this is clear.
If $x_3 \in [0,1]$, it relies on the definition of $a$ 
in the last inequality of the chain 
\begin{align*}
x_1^2 x_2^2 x_3 & \le \f{(x_1^2+ x_2^2)^2}{4} x_3 
\le \f{(1- x_3^2)^2}{4} x_3
= (1-x_3) \bigg[ \f{(1+x_3)^2(1-x_3)x_3}{4} \bigg]\\
& \le (1-x_3) \Big[ a (1 + 4 x_3 + 4 x_3^2) \Big]
= a (1-T_3(x_3)).
\end{align*}
Altogether, we have obtained $| (m_{(2,2,1)} - p^*)(x)| = \| m_{(2,2,1)} - p^* \|_B = a$ for all $x \in \cS$, as announced.

The second step comprises showing that $\| m_{(2,2,1)} - p \|_B \ge a$
for any $p \in \cP_4^3$,
or in fact for some best approximant $p$ to $m_{(2,2,1)}$ from $\cP_4^3$.
Let us momentarily take for granted that $v = p^* - p$ satisfies $v(x) \ge 0$ for some $x \in \cS^+$.
With the help of this $x \in \cS^+$, we derive
$$
\|m_{(2,2,1)} - p \|_B
\ge  (m_{(2,2,1)} - p)(x)
=  (m_{(2,2,1)} - p^*)(x) + (p^* - p)(x)
=  a + v(x) \ge a,
$$
as desired.
Thus, it remains to justify that existence of $x \in \cS^+$ such that $v(x) \ge 0$.
According to Proposition \ref{PropRedu2},
the best approximant $p$ can be chosen to inherit properties of $m_{(2,2,1)}$,
in particular being odd in $x_3$,
even in $x_1$ and $x_2$,
and symmetric when swapping $x_1$ and $x_2$.
Therefore, one can choose $p$ to contain only the terms $x_3$, $x_3^3$, and $(x_1^2+x_2^2)x_3$,
and when restricted to the boundary of~$B$,
it contains only the terms $x_3$ and $x_3^3$, just like $p^*$.
As a consequence, we write, for some $c,d \in \bR$,
$v(x) = c x_3 + d x_3^3$ for all $x \in \cS^+$.
Now assume by contradiction that $v(x) < 0$ for all $x \in \cS^+$.
This translates, after simplification, into
$c+d<0$, $-4c - d < 0$, and $\tau^{-2}c + d<0$.
Adding the second to the first yields $-3c<0$,
while adding the second to the third yields $(\tau^{-2}-4)c<0$,
which is impossible since $\tau^{-2} \approx 6.088 > 4$.
Thanks to this contradiction, the proof is complete.
\epf

We have seen that the signature in the above proof was supported on the boundary of the ball.
As a matter of fact,
this is a phenomenon we noticed in all cases (where we could extract signatures).
We thus conjecture that,
for the approximation of monomials on the euclidean ball,
signatures actually live on the sphere.

\subsection{The cross-polytope}

To the best of our knowledge,
Chebyshev polynomials relative to the cross-polytope have not been investigated in the literature,
hence the values of $E(k,C)$ reported in Table \ref{TableCrossPolytope} seem to all be new.
As expected, they are smaller than the values of $E(k,B)$ presented in Table \ref{TableBall} and larger than the values of $E(k,S)$ shown in Table \ref{TableSimplex}.

\begin{table}[h]
\small
\begin{center}
\begin{tabular}{|c|c|c|c|}
\hline
degree $n=3$ ($\times 10^{-2}$) & degree $n=4$ ($\times 10^{-2}$) & degree $n=5$ ($\times 10^{-3}$) & degree $n=6$ ($\times 10^{-3}$)\\
\hline
\hline
$\mathbf{E((1,1,1),C) \approx 3.703}$ & $\mathbf{E((2,1,1),C) \approx 1.273}$ & $\mathbf{E((3,1,1),C) \approx  4.764}$ & $\mathbf{E((4,1,1),C) \approx 1.853}^*$ \\
\hline
 & & $\mathbf{E((2,2,1),C) \approx 3.398}$ & $\mathbf{E((3,2,1),C) \approx 1.087}\phantom{^*}$\\
\hline
 & & & $\mathbf{E((2,2,2),C) \approx 0.661}^*$\\
 \hline 
\end{tabular}
\end{center}
\caption{Cross-polytope in dimension $d=3$:
all values were previously unknown.}
\label{TableCrossPolytope}
\end{table}

\subsection{The simplex}

When $\Omega = S$,
the values of $E(k,\Omega)$ had earlier been found when $|k|=3$, 
but not for $|k|=4$, $|k|=5$, and $|k|=6$ (except $k=(2,2,2)$).
All these values are shown in Table \ref{TableSimplex}.
It appears empirically that signatures (when they could be extracted) live on the boundary of the domain,
and more precisely here on the face of equation $x_1+x_2+x_3 = 1$.
Note that this would imply,
due to the close connection between the approximation of $m_k$ on $S$ and the approximation of $m_{2k}$ on~$B$, 
some particular cases of the conjecture relative to $B$.
    
\begin{table}[h]
\small
\begin{center}
\begin{tabular}{|c|c|c|c|}
\hline
degree $n=3$ ($\times 10^{-2}$) & degree $n=4$ ($\times 10^{-3}$) & degree $n=5$ ($\times 10^{-4}$) & degree $n=6$ ($\times 10^{-4}$)\\
\hline
\hline
$E((1,1,1),S) \approx 1.388$ & $\mathbf{E((2,1,1),S) \approx 2.688}$ & $\mathbf{E((3,1,1),S) \approx 5.984^*}$ & $\mathbf{E((4,1,1),S) \approx 1.405^*}$ \\
\hline
 & & $\mathbf{E((2,2,1),S) \approx 4.695}\phantom{^*}$ & $\mathbf{E((3,2,1),S) \approx 1.000^*}$\\
\hline
 & & & $E((2,2,2),S) \approx 0.6265\phantom{^*}$\\
 \hline 
\end{tabular}
\end{center}
\caption{Simplex in dimension $d=3$:
the previously unknown values are shown in boldface.}
\label{TableSimplex}
\end{table}

In the case $k=(2,1,1)$,
we could use the insight brought forward by our computations
to derive a Chebyshev polynomial explicitly.
The result reads as follows.

\bthm
\label{ThmS211}
With $\tau \in [0,1/4]$ being the solution to $\max \big\{ y(1-2y)(y-\tau)^2, y \in [0,1/2] \big\} = \tau^2/18$ and with $c \coloneqq -3/\tau$,
the error of best approximation on the simplex to $m_{(2,1,1)}(x_1,x_2,x_3) = x_1^2 x_2 x_3$ by trivariate polynomials of degree at most $3$ is 
$$
E_{\cP_3^3}(m_{(2,1,1)}, S) = \f{1}{2 c^2},
$$
while a Chebyshev polynomial is given by 
\begin{align*}
P(x) & = 
x_1^2 x_2 x_3 \\
& + \f{1}{2 c^2}
\left[ -16 x_1^2 (x_2 + x_3) + 16 x_1 (x_2+x_3)^2 - 2 (64 + 12 c + c^2) x_1 x_2 x_3 +  8 x_2 x_3 - 2 (x_2 + x_3)  + 1 \right].
\end{align*}
\ethm

\bpf
Guided by our computational procedure,
we make the guess---which we are about to verify---that there is a signature with support
$\cS = \cS^+ \cup \cS^-$,
where 
\begin{align*}
\cS^+ & = \left\{
\bpmx 1/4 \\ 3/4 \\ 0 \epmx,
\bpmx 1/4 \\ 0 \\ 3/4 \epmx,
\bpmx 0 \\ 1/2 \\ 1/2 \epmx,
\bpmx 1 \\ 0 \\ 0 \epmx,
\bpmx 1-2\tau \\ \tau \\ \tau \epmx
\right\},\\
\cS^- & = \left\{
\bpmx 3/4 \\ 1/4 \\ 0 \epmx,
\bpmx 3/4 \\ 0 \\ 1/4 \epmx,
\bpmx 0 \\ 0 \\ 1 \epmx,
\bpmx 0 \\ 1 \\ 0 \epmx,
\bpmx 1-2\sigma \\ \sigma \\ \sigma \epmx
\right\}.
\end{align*}
The parameters $\tau,\sigma \in [0,1/2]$ are considered free for now---their specific choice will be revealed later.

The first step consists in proving that $|P(x)| = \|P\|_S$ for all $x \in \cS$.
To this end, we start by proving that $\|P\|_S \le 1/(2 c^2)$,
i.e., that $|P(x_1,x_2,x_3)| \le 1/(2c^2)$ whenever $x_1,x_2,x_3 \ge 0$ and $x_1+x_2+x_3 \le 1$.
This is easy to see if $x_1 = 0$, since then $2 c^2 \,P(0,x_2,x_3) 
= 8 x_2 x_3 - 2 (x_2 + x_3)  + 1$
is bounded as
$$
\left\{
\bmx
\ge -2 (x_2+x_3) + 1 \ge -2 + 1 \phantom{+\,} &= -1,\\
\le 2 (x_2 + x_3)^2 - 2 (x_2 + x_3) + 1 &\le 1.
\phantom{-}
\emx
\right.
$$
Using the fact that $8t(1-2t) \le 1$ for all $t \in \bR$,
it is also easy to see that $|P(x_1,x_2,x_3)| \le 1/(2c^2)$ if $x_2=0$ or $x_3=0$,
since, e.g., $2 c^2 \,P(x_1,x_2,0) = -16x_1^2x_2 + 16 x_1 x_2^2 - 2 x_2 +1$ is
$$
= \left\{
\bmx
-2 x_1 [ 8x_2(x_1-x_2)] -2 x_2 + 1
&\ge& -2 x_1 [ 8x_2(1-2x_2)] -2 x_2 + 1
\ge -2x_1 - 2x_2 + 1 \ge -1 ,\\
\phantom{-}
2 x_2 [ 8x_1(x_2-x_1)] -2 x_2 + 1
&\le& 
\phantom{-\,}
2 x_2 [8x_1(1-2 x_1)] - 2 x_2 + 1
\le \phantom{-\,} 2x_2 - 2x_2 + 1
= 1. \phantom{-\,}
\emx
\right.
$$
We therefore consider the case where $x_1,x_2,x_3$ are all nonzero 
and we notice that we can assume $x_2=x_3$.
Indeed, if we reparametrize by setting $y = (x_2+x_3)/2$
and $z = (x_2-x_3)/2$,
so that $x_2x_3 = y^2 - z^2$,
then the expression for $P(x)$ becomes
\begin{align}
\label{ExprPxSym}
P(x) & = x_1^2 y^2 + \f{1}{2c^2}  \left[ -32 x_1^2 y + 64 x_1 y^2 - 2 (64 + 12 c + c^2) x_1 y^2 +  8 y^2 - 4 y  + 1 \right]\\
\nonumber
& - x_1^2 z^2+ \f{1}{2c^2} \left[ 
\phantom{ -32 x_1^2 y + 64 x_1 y^2 -\;\;}
2 (64 + 12 c + c^2) x_1 z^2 -  8 z^2 
\phantom{\;\;- 4 y  + 1}\right]\\
\nonumber
& = \f{1}{2c^2} \left[
2 c^2 x_1^2 y^2 -32 x_1^2 y - 2 (32 + 12 c + c^2) x_1 y^2 +  8 y^2 - 4 y  + 1
\right]
- z^2 q(x_1)
\end{align}
for some univariate quadratic polynomial $q$.  
Thus, given $x \in S$ such that $|P(x) | = \|P\|_S$ and $t \in \bR$ small enough,
we define $x^{(t)} \in S$  by
$x^{(t)}_1 = x_1$,
$x^{(t)}_2 = x_2+t$,
and $x^{(t)}_3 = x_3-t$.
In view of $(x^{(t)}_2+x^{(t)}_3) = (x_2+x_3)/2 \eqqcolon y$
and of $(x^{(t)}_2-x^{(t)}_3) = (x_2-x_3)/2 + t \eqqcolon z + t$, 
while supposing e.g. that $P(x)>0$, 
the inequality $P(x^{(t)}) \le P(x)$
reads
$P(x_1,y,y) - (z+t)^2 q(x_1) \le
P(x_1,y,y) - z^2 q(x_1)$ whenever $|t|$ is small enough.
This implies that $z q(x_1) = 0$,
hence  that $P(x) = P(x_1,y,y)$,
meaning that the last two coordinates of an extremal point can be assumed to be equal, as claimed.
Now, to determine the maximum of $|P(x_1,y,y)|$ when $x_1,y \ge 0$ and $x_1+2y \le 1$,
we recall from \eqref{ExprPxSym} that
\be
\label{ExprPxSym2}
P(x_1,y,y) = \f{1}{2c^2} \left[
2 c^2 x_1^2 y^2 -32 x_1^2 y - 2 (32 + 12 c + c^2) x_1 y^2 +  8 y^2 - 4 y  + 1
\right],
\ee
so that 
$$
\f{\partial P(x_1,y,y)}{\partial y}
= 
\f{1}{2c^2} \left[
4 c^2 x_1^2 y -32 x_1^2 - 4 (32 + 12 c + c^2) x_1 y +  16 y - 4 
\right] .
$$
As a consequence, at a critical point, we have 
$(32 + 12 c + c^2) x_1 y =  c^2 x_1^2 y - 8 x_1^2  +  4 y - 1 $ and in turn
$$
P(x_1,y,y) = 
\f{1}{2c^2} \left[
-16 x_1^2 y - 2 y +1
\right].
$$
It follows that $|P(x_1,y,y)| \le 1/(2c^2)$ at any critical point,
since $[-16 x_1^2 y - 2 y +1 ] \le 1$ is clear,
while $[-16 x_1^2 y - 2 y +1 ] \ge -1$ holds because $16 x_1^2 y + 2 y = 2y(8 x_1^2 + 1) \le (1-x_1)(8 x_1^2 + 1)$,
the latter having maximal value $(68+5\sqrt{10})/54 \approx 1.5520 \le 2$ over $x_1 \in [0,1]$.
At this point, it remains to verify that $|P(x_1,y,y)| \le 1 / (2c^2)$
on the boundary of the domain $\{ x_1 \ge 0, y \ge 0, x_1 + 2y \le 1\}$.
Since the cases $x_1=0$ and $y=0$ have already been dealt with,
we need to consider the case $x_1 = 1-2y$, $y \in [0,1/2]$.
After some technical calculations left to the reader,
starting from \eqref{ExprPxSym2} and recalling that $\tau = - 3/c$,
we arrive at 
\be
\label{ExprPxSym3}
P(1-2y,y,y)  = \f{1}{2 c^2} - 2 y (1-2y)(y-\tau)^2, 
\qquad  y\in [0,1/2].
\ee
The inequality $P(1-2y,y,y) \le 1 / (2 c^2)$ is obvious from here.
The parameter $\tau$ is chosen to secure the other inequality $P(1-2y,y,y) \ge -1 / (2 c^2)$,
which is equivalent to $y (1-2y)(y-\tau)^2 \le 1/(2 c^2)$ for all $y \in [0,1/2]$
and thus follows from the equation $\max \big\{ y(1-2y)(y-\tau)^2, y \in [0,1/2] \big\} = \tau^2/18$.
Note that this equation has a (unique) solution in $[0,1/4]$,
because $\tau^2/18$ increases from $0$ to $1/288$ on this interval and  $\max \big\{ y(1-2y)(y-\tau)^2, y \in [0,1/2] \big\}$ decreases from a positive quantity to $1/512$ there.
In consequence,
we have now established that $|P(x)| \le 1/(2c^2)$ for all $x \in S$, as announced.\\
Let us turn to the justification that $|P(x)| = 1/(2c^2)$ for all signature points $x \in \cS = \cS^+ \cup \cS^-$.
From \eqref{ExprPxSym3},
we immediately see that $P(x) = 1/(2c^2)$ for the last three points of $\cS^+$.
The fact that $P(x) = -1/(2c^2)$ for the last point of $\cS^-$ is due to the choice of $\sigma$,
as we take it to be the maximizer of $y(1-2y)(y-\tau)^2$ over $y \in [0,1/2]$,
ensuring that $\sigma (1-2\sigma)(\sigma-\tau)^2 = 1/(2c^2)$ and hence that $P(1-2\sigma,\sigma,\sigma) = -1/(2 c^2)$.
As for the other signature points,
notice that they are of the form $(1-y,y,0)$ or $(1-y,0,y)$,
for which technical calculations left to the reader yield
$$
P(1-y,y,0)  = P(1-y,0,y) = - \f{1}{2 c^2} T_3(2y-1),
\qquad y \in [0,1],
$$
where, as usual, $T_3(t) = 4 t^3 - 3t$ is the univariate Chebyshev polynomial of degree $3$.
It then easily follows that $P(x) = 1/(2c^2)$ for the first two points of $\cS^+$
and that $P(x) = -1/(2c^2)$ for the first four points of $\cS^-$.   
Altogether,
we have now obtained $|P(x)| = \|P\|_S = 1/(2c^2)$ for all $x \in \cS$, as announced.

The second step consists in proving that $\|m_{(2,1,1)} - p\|_S \ge 1/(2c ^2)$ for any $p \in \cP_3^3$.
To this end, we notice that all signature points $x \in \cS$ lie on the face $\cF \coloneqq \{x_1 \ge 0, x_2 \ge 0, x_3 \ge 0, x_1 + x_2 + x_3 =1 \}$
and we remark that $\{ (x_1,x_2): x \in \cS \}$
is the support of an extremal signature for $\cP_3^2$ associated with $P_{| \cF}$,
in the sense that there exist $c_x > 0$, $x \in \cS$,
such that $\sum_{x \in \cS} c_x \, \sgn(P(x)) \, r(x_1,x_2) = 0$ for all $r \in \cP_3^2$.
This can be verified (numerically) by looking at the null space of the $10 \times 10$ matrix with entry $\sgn(P(x)) m_{(k_1,k_2)}(x_1,x_2)$
on the row indexed by $(k_1,k_2)$ with $k_1+k_2 \le 3$ and on the column indexed by $x \in \cS$.
Theorem \ref{ThmSig} could now be invoked.
Alternatively,
given $p \in \cP_3^3$,
we can write $m_{(2,1,1)} - p = P - q$ for some $q \in \cP_3^3$
and consider an $x \in \cS$ such that $\sgn(P(x)) \, q(x) \le 0$,
which is possible for otherwise $\sum_{x \in \cS} c_x \, \sgn(P(x)) \, q(x_1,x_2,1-x_2-x_3) = 0$ would be violated.
In this way, the desired inequality follows from
$$
\|m_{(2,1,1)} - p\|_S \ge |P(x) - q(x)| 
\ge |P(x)| = \f{1}{2c^2}.
$$
The proof is now complete.
\epf

\brk
The numerical values of the parameters $\tau$ and $\sigma$ are  $\tau \approx 0.21998$ and $\sigma \approx 0.41942$,
leading to the error of best approximation the numerical value $E_{\cP_3^3}(m_{(2,1,1)}, S) \approx 2.68850\times 10^{-3}$. 
In fact,
it can be shown (a computer algebra system will facilitate the task)
that $\tau$ is the smallest real root of the quartic polynomial $2880 t^4 - 5472 t^3 + 4880 t^2 -1944 t + 243$.
\erk

\section{Conclusion}
\label{SecConc}

In this article,
we proposed a semidefinite-programming method to compute best approximants to monomials by multivariate polynomials of lower degree.
More than providing numerical values,
the method  allows us to make guesses for the multivariate analogs of Chebyshev polynomials
that can---sometimes---be {\em a posteriori} certified explicitly or symbolically.
Of note, the generic nonuniqueness of such analogs prompted us to preferentially solve the dual optimization program,
putting the classical notion of signature at center stage.
We emphasize that the underlying methodology is quite versatile and should enable to attack other problems in multivariate Approximation Theory by relying on modern tools from Optimization Theory,
so long as one is ready to give up on purely analytical solutions.
This is a point of view already brought forward by a subset of the authors to determine minimal projections (exploiting moments, see \cite{FouLas1}) and Chebyshev polynomials associated to union of intervals (exploiting sums-of-squares, see \cite{FouLas2}).
In the multivariate setting,
we should also be able to deal with unions of domains,
as well as tackling norms other than $L_\infty$ 
(notably $L_1$ and $L_{2m}$, $m \in \bN$),
adding convex constraints (e.g. interpolatory, shape-enforcing, etc),
in the spirit of the proof-of-concept software {\sf Basc},
short for `Best Approximations by Splines under Constraints', see \cite{Basc}.
We note, though, that the semidefinite programs encountered in {\sf Basc} could only be handled 
thanks to the benefits of representing univariate polynomials in the Chebyshev basis rather than in the monomial basis,
so a similar approach should be taken in the multivariate setting.
This is indeed realizable, at least in theory,
and we give pointers on how to do this in some supplementary material.
Still, bringing {\sf Basc} to the multivariate realm will be a mighty task, but certainly one worth taking by a fresh generation of approximators/optimizers.

\paragraph{Acknowledgment.}
This work is the result of a collaboration made possible by the SQuaRE program at the American Institute of Mathematics (AIM).
We are truly grateful to AIM for 
the supportive and mathematically rich environment they provided.
In addition, we owe thanks to several funding agencies,
as 
M. D. is supported by the Australian Research Council Discovery Early Career Award DE240100674, 
S. F. is partially supported by grants from the National Science Foundation (DMS-2053172) and from the Office of Naval Research (N00014-20-1-2787),  
E. de K. is supported by grants from the Dutch National Science Foundation (NWO) 
(OCENW.M.23.050 and OCENW.GROOT.2019.015),
J. B. L. is supported by the AI Interdisciplinary Institute  through the French program ``Investing for the Future PI3A" (ANR-19-PI3A-0004)
and by the National Research Foundation, Singapore,
through the DesCartes and Campus for Research Excellence and Technological Enterprise (CREATE) programs, 
and Y. X. is partially supported by the Simons Foundation (grant \#849676).


\begin{thebibliography}{99}

\bibitem{AV1} 
N. N. Andreev and V. A. Yudin.
{\em Polynomials of least deviation from zero and Chebyshev-type 
cubature formulas.} 
Proceedings of the Steklov Institute of Mathematics,  232,  39--51, 2001.

\bibitem{AV2} 
N. N. Andreev and V. A. Yudin.
{\em Best approximation of polynomials on the sphere and on the ball.}
In:  Recent Progress in Multivariate Approximation,
International Series of Numerical Mathematics  137.
Birkh\"auser, 2001.

\bibitem{BHN}
 B. D. Bojanov, W. Haussmann, and G. P. Nikolov.
{\em  Bivariate polynomials of least deviation from zero.}
Canadian Journal of Mathematics,  53, 489--505,  2001.

\bibitem{CheOri}
P. L. Chebyshev.
{\em Sur les questions de minima qui se rattachent \`a la repr\'esentation approximative des fonctions.}
M\'emoires de l'Acad\'emie Imp\'eriale des Sciences de St.-Petersbourg, 7, 199--291, 1859. 

\rev{
\bibitem{DFKJLX}
M. Dressler, S. Foucart, M.  Joldes, E. de Klerk, J.~B. Lasserre, and Y. Xu.
{\em Least multivariate Chebyshev polynomials on diagonally determined domains.}
 arXiv:2405.19219v1 [math.OC].
}
\bibitem{EhlZel}
H. Ehlich and K. Zeller. 
{\em $\check{\mbox{C}}$eby$\check{\mbox{s}}$ev-Polynome in mehreren Ver\"{a}nderlichen.} 
Mathematische Zeitschrift, 93, 142--143, 1966.

\bibitem{FouLas1}
S. Foucart and J. B. Lasserre. 
{\em Determining projection constants of univariate polynomial spaces.} 
Journal of Approximation Theory, 235, 74--91, 2018.

\bibitem{FouLas2}
S. Foucart and J. B. Lasserre. 
{\em Computation of Chebyshev polynomials for union of intervals.} 
Computational Methods and Function Theory, 19/4, 625--641, 2019.

\bibitem{Basc}
S. Foucart and V. Powers. 
{\em {\sf Basc}: constrained approximation by semidefinite programming.} 
IMA Journal of Numerical Analysis, 37/2, 1066--1085, 2017.

\bibitem{G}
W. B. Gearhart.
{\em Some Chebyshev approximations by polynomials in two variables.}        
Journal of Approximation Theory,  8, 195--209, 1973.

\revv{
\bibitem{GloptiPoly solution extraction}
D. Henrion and J.B. Lasserre. 
{\em Detecting global optimality and extracting solutions in GloptiPoly.} 
In: Positive Polynomials in Control,
Lecture Notes in Control and Information Science 312.
Springer, 2005.
}

\bibitem{GloptiPoly}
D. Henrion, J. B. Lasserre, and J. Loefberg. 
{\em GloptiPoly 3: moments, optimization and semidefinite programming. }
Optimization Methods and Software,  24/4-5,  761--779, 2009.
\url{https://homepages.laas.fr/henrion/software/gloptipoly3/}

\revv{
\bibitem{lass-siopt-01}
J.B. Lasserre.
{\em Global optimization with polynomials and the problem of moments}. {SIAM Journal on Optimization} 11, 796--817, 2001.
}

\bibitem{JBL}
J. B. Lasserre. 
{\em The Moment-SOS Hierarchy.} 
In: Proceedings of the International Congress of Mathematicians (ICM 2018). 
World Scientific, pp. 3773--3794, 2019.

\bibitem{LauRen}
\rev{
M. Laurent and F. Rendl. 
{\em Semidefinite Programming  and Integer Programming}.
In: Handbooks in Operations Research and Management Science.
Elsevier, pp. 393--514, 2005.}

\bibitem{MoaPeh}
I. Moale and F. Peherstorfer. 
{\em An explicit class of min–max polynomials on the ball and on the sphere.} 
Journal of Approximation Theory 163/6, 724--737, 2011.

\bibitem{NewXu}
D. J. Newman and Y. Xu. 
{\em Tchebycheff polynomials on a triangular region.} 
Constructive Approximation, 9/4, 543--546, 1993.

\rev{
\bibitem{Put}
M. Putinar. 
{\em Positive polynomials on compact semi-algebraic sets.} 
Indiana University Mathematics Journal, 42/3, 969--984, 1993.
}

\bibitem{R}
 M. Reimer.
{\em On multivariate polynomials of least deviation from zero on the unit ball.}
Mathematische Zeitschrift,  153, 51--58, 1977. 

\bibitem{Rich}
H. Richter. \emph{Parameterfreie Absch\"atzung und Realisierung von Erwartungswerten}, Deutsche Gesellschaft f\"ur Versicherungsmathematik, 3, 147--162, 1957.

\bibitem{Riv}
T. J.  Rivlin. 
Chebyshev Polynomials. 
Second edition. 
Courier Dover Publications, 2020.

\bibitem{shapiro}
T. J. Rivlin and H. S. Shapiro.
{\em A unified approach to certain problems of approximation and optimization.}
SIAM Journal on Numerical Analysis, 9,  670--699, 1961.

\rev{
\bibitem{Schm}
K. Schm\"udgen. 
The Moment Problem. Springer, 2017.
}

\bibitem{Slo}
J. Sloss. 
{\em Chebyshev approximation to zero.} 
Pacific Journal of Mathematics, 15/1, 305--313, 1965.


\bibitem{VY}
V. A. Yudin, 
{\em Best approximation to monomials on a cube.} Sbornik Mathematics, 199, 1251--1262, 2008.

\bibitem{X04}
Y. Xu.
{\em On polynomials of least deviation from zero in several variables.}
\rev{Experimental Mathematics}, 13, 103--112, 2004. 

\bibitem{X05}
Y. Xu.
{\em Best approximation of monomials in several variables.}
Rendiconti del Circolo Matematico di Palermo Series 2, 76, 129--155,  2005. 

\end{thebibliography}
\end{document}